\documentclass[a4paper]{article}
\usepackage{graphicx}
\usepackage{amssymb}

\def\RR{{\mathbb R}}

\def\CC{{\mathbb C}}
\def\HH{{\mathbb H}}

\def\PP{{\mathbb P}}
\def\bra{{\langle}}
\def\ket{{\rangle}}

\def\be{\begin{equation}}
\def\ee{\end{equation}}
\def\reel{{\rm Re}}

\newcommand{\bigL}{\mathop{\mathrm{Lex}}}

\newtheorem{propo}{Proposition}

\newtheorem{definition}{Definition}

\newtheorem{theorem}{Theorem}

\begin{document}
\title{Two roads to noncommutative causality}

\author{Fabien Besnard\footnote{Pôle de recherche M.L. Paris, EPF, 3 bis rue Lakanal, F-92330 Sceaux. \goodbreak fabien.besnard@epf.fr}}

\maketitle

\begin{abstract}
We review the physical motivations and the mathematical results obtained so far in the isocone-based approach to noncommutative causality. We also give a briefer account of the alternative framework of Franco and Eckstein which is based on Lorentzian spectral triples. We compare the two theories on the simple example of the product geometry of the Minkowski plane by the finite noncommutative space with algebra $M_2(\CC)$.
\end{abstract}

\section{Introduction}

The properties of the electromagnetic field led the early 20th century physicists to a dramatic change in their notions of space and time. Yet more upheaval came from Einstein's theory of the gravitational field. Today, known fields other than gravitation are described together in the Standard Model of particle physics. The derivation of this very complex piece of machinery from noncommutative geometric ingredients is a major breakthrough \cite{cm}. While it sheds light on some of the most contrived aspects of the Standard Model, such as the Higgs mechanism, it also opens the door  to a renewed understanding of the structure of spacetime. However two hard problems remain. The first is the passage from the euclidean signature  to the physical lorentzian one, which is still elusive for noncommutative manifolds (in spite of recent progress in some particular cases, \cite{vddg}). The second is the quantization of the noncommutative scheme \cite{besq}.

In fact, there is a third difficulty, which we believe is the key to the elucidation of the first  two:  the physical interpretation of the theory. In noncommutative geometry, a $C^*$-algebra is interpreted as an algebra of would-be functions on some virtual space. In the case of the spectral Standard Model the algebra ${\cal A}$ is composed of true  functions on space with values in $\CC\oplus\HH\oplus M_3(\CC)$.  On a different stance, coordinates are expected to become noncommutative observables in quantum gravity \cite{dfr}. We think that it is a fruitful point of view to consider the elements of ${\cal A}$ already as quantum observables, which take into account the  quantum nature of the ``inner'' space unraveled by the Standard Model, but not yet the quantum nature of the spacetime manifold.

With this in mind, we try a causal approach to the problem of the lorentzian signature. More precisely, we observe that causality is a   property which is specific to lorentzian manifolds, and that moreover the causal order relation encapsulates the metric information up to a conformal factor. This is hardly a new observation: it is a the root of the Causal Sets approach to quantum gravity \cite{sorkin}. Thus, we want to characterize, among the general elements of the algebra of noncommutative  functions ${\cal A}$, those which play the role of causal functions. The subset of these distinguished elements is called an isocone. The isocone program has been started in \cite{bes1} and is still in its early stages of developments. In particular it does not yet make contact with the Dirac operator, which is the central tool of   noncommutative geometry. In other words it is not yet known how to distinguish in the noncommutative case true causality, which comes from a metric, from general order relations. 

However there is another theory, started by Franco and Eckstein in \cite{fe1}, which aims at defining causality out of a Dirac operator. Of course this top-down approach depends on a correct formulation of noncommutative geometry in the lorentzian signature in the first place.

In this paper we will review the isocone-based approach to causality, starting with general physical motivations in section 2. This will lead us to a set of four fundamental hypotheses, which we will transform in section 3  into a   mathematically more tractable set of axioms that  defines isocones in the $C^*$-algebra category. In section 4 we will review  general properties of isocones, and in section 5 some classification results for them. In section 6 we briefly review the Franco-Eckstein theory, and compare the two approaches on a simple example in section 7. We will then conclude with a few remarks.

In all the text ``compact'' means ''compact  and Hausdorff'' and the order relation used on $\RR$ is always the natural one.

    
\section{Time, clocks and causality} 

We start with a rather philosophical discussion, where ``philosophical'' is here an elegant substitute for ``vague''. This vagueness will perhaps be forgiven since the question we want,  to address ``what is time ?'',  is certainly one of the most difficult of all. It seems reasonable to look for an answer in our two best theories of natural phaenomena to date, namely General Relativity and Quantum theory. The former theory carefully distinguishes between chronology, duration and causality, and puts strong structural constraints on them. The latter theory says nothing at all: time in quantum theory is an external parameter. So the punchline appears to be this: time, in any of its three facets mentioned above, is structurally defined by GR and has no quantum aspect at all. Note that by ``structurally defined'' we mean that there is no more explicit definition of time in GR than there is any explicit definition of sets in axiomatic set theory, but that the two theories implicitly define their objects by specifying how they behave. That the definition of time is a structural one is not really surprising, the only alternative being that it is an emergent phenomenon that can be defined in term of something else. However it is unsatisfying that quantum theory should have no word on this matter, and it is rather unlikely that we have here the definitive answer.

There are indeed several arguments in favor of an influence of quantum theory on the question  of time, and we refer to \cite{tppm} for some of them. Here we quote only two. The first and most obvious is that if the gravitational field has to be quantized, then the structural constraints on time observables will have to be modified accordingly. The second is that the so-called clock hypothesis, which is a postulate of GR, states that a clock measures the proper time of an observer who carries it. But clocks being made of quantum fields, like everything else, cannot in principle be properly modelized without the help of quantum theory. In other words there should be quantum observables associated to them. A possible objection to this point of view would be that clocks are necessarily macroscopic objects. Once again this would amount to say that time is an emergent phenomenon.

The characterization of quantum observables associated to measures of duration poses a real challenge, even in the context of quantum mechanics on a non-relativistic background \cite{tqm}. However, we will argue that if we fall back to the cruder notion of causality, a tentative answer to this problem can be given. Let us now explain more precisely what we propose to do.

From General Relativity we borrow the concept of \emph{causality}. A spacetime $M$, i.e. a connected time-orientable Lorentzian manifold, is said to be \emph{causal} if the following relation

\be
p\preceq q \Leftrightarrow \mbox{ there exists a causal curve from }p\mbox{ to }q\label{causalorder}
\ee  

is a partial order. Not all metrics define causal spacetimes, even those which solve the Einstein's equations can have pathological causal properties. However, there always exists a neighbourhood $U$ of a given event $p$ such that $\preceq$ is a partial order on $U$. In fact much more is true: one can find a neighbourhood with the best possible causality property, namely global hyperbolicity (see \cite{mingsan}, th 2.14).

The causal ordering permits to define \emph{causal functions}, that is functions $f : M\rightarrow \RR$ such that $x\preceq y\Rightarrow f(x)\le f(y)$. These functions are classical observables which correspond to measurement devices that we call \emph{generalized clocks}. A generalized clock, if compared to a genuine clock, can arbitrarily change pace or even stop for a while. The only thing it cannot do is run backwards. We now want to associate quantum observables to these generalized clocks.  For this we introduce a quite general setting for quantum physics \cite{emch}.

\begin{enumerate}
\item For a given physical system, there exists a set of states ${\cal S}$, a set of observables ${\cal O}$, and a  duality map  ${\cal S}\times {\cal O}\rightarrow \RR$, which is separating for ${\cal S}$ and ${\cal O}$.
\item ${\cal S}$ is a convex set.
\item\label{somme} ${\cal O}$ is a real vector space.
\item\label{VNrule} ${\cal O}$ is stable under $a\mapsto f(a)$ where $f$ is a real-valued function.
\end{enumerate}

Note that we do not assume for the moment that the duality function gives the expectation value of the observable when the system is in some state: this is a matter of interpretation of the formalism and we will deal with that later on. There are many good reasons to believe in these axioms, and we refer to \cite{emch} for justifications. Let us only comment on the last  of them, which will be crucial for our purpose. This axiom follows from the following intuitive interpretation: the observable $f(a)$ corresponds to an observation of $a$ followed by the mathematical operation of applying the function $f$ to any result of such an observation.

A neat way of fulfilling all these requirements is to suppose that ${\cal O}$ is the self-adjoint part of a $C^*$-algebra and that ${\cal S}$ is its state space, this is the so-called $C^*$-paradigm. However it is not quite the most general solution: a Jordan algebra would suffice. So we will start by assuming that ${\cal O}$ is a Jordan algebra. We postpone to the next section the specification of the type of Jordan algebra that we will use in order to make sense of our hypotheses.



 We now call ${\cal C}$ the set of \emph{causal observables} corresponding to generalized clocks. Our first hypothesis, justified above, is that ${\cal C}$ is a subset of a Jordan algebra.
 
\begin{quote}
{\bf Hypothesis 0 (Jordan paradigm)} ${\cal C}$ is a subset of a Jordan algebra ${\cal O}$.
\end{quote}

In some classical limit ${\cal C}$ should correspond to the set of causal functions on $M$,  but it should be stressed that the correspondence is in the direction causal observable $\rightarrow$ causal functions, in particular we do not assume the existence of a quantization map which would go in the other direction. In fact it is well-known that quantization maps cannot commute with functional calculus and be linear at the same time \cite{ae}. What are the properties of ${\cal C}$ ? Consider $o\in {\cal C}$, $c$ a generalized clock corresponding to $o$ and $f$ a non-decreasing real-valued function. Then the device $c'$ which is formed by $c$ followed by a computation of the function $f$ on the result obtained by $c$ is still a generalized clock. Hence $f(o)\in {\cal C}$. This is our hypothesis number 1 for the set ${\cal C}$.
 
\begin{quote} {\bf Hypothesis 1 (functional stability)} If $o\in{\cal C}$ and $f$ is a real non-decreasing  function, then $f(o)\in{\cal C}$. 
\end{quote}

Let us remark here that this implies that the algebra is unital and that ${\cal C}$ contains all the real multiples of the unit. Let us dissipate any worry that could arise at this point: it does not mean that our approach only applies to compact spacetimes. If we start with some noncompact spacetime and use the approach advocated here with a unital commutative algebra, we will recover some compactification of the spacetime we started from, depending on the algebra used. Of course   these objects have boundaries. It is indeed a well-known theorem that compact causal spacetimes always have boundaries (or to put it another way, that compact closed causal spacetimes do not exist).


We now introduce a separation hypothesis.

\begin{quote} {\bf Hypothesis 2 (separation)} ${\cal C}$ separates the pure states of ${\cal O}$.
\end{quote}
 
The origin of this hypothesis is the classical situation. In any sufficiently small given region of spacetime one can introduce causal coordinates, for instance GPS coordinates \cite{rov}, which obviously separates events (i.e. pure states). If one considers the spacetime manifolds globally, there can be too few causal functions to separate the events (in fact, there can be none at all, except the constants). However the hypothesis will be fulfilled if one assumes strong causality, which despite its name, is a weaker causality condition than global hyperbolicity \cite{mingsan}.


Hypothesis 2 can thus be viewed as an exportation to the quantum regime of the idea that causality plays a fundamental role in the process of identification of what we call ``different spacetime events''. As such it cannot be fully justified, it is an act of faith.

As soon as we have adopted it, this hypothesis grants us a partial ordering on the set $P({\cal O})$ of pure states of ${\cal O}$, which is so defined:

\be
\forall p_1,p_2\in P({\cal O}), p_1\preceq p_2\Leftrightarrow p_1(c)\le p_2(c)\mbox{ for all }c\in{\cal C}\label{deforder}
\ee

In the classical situation, if we started from the set of causal functions, we recover in this way the causal order on  spacetime events. Mind that this is not a mere tautology: it could well be that the order defined through (\ref{deforder}) turn out to be stronger than the causal order (see the definition of toposets in section 4). But as is now familiar, the equality between causal order and the order defined through the causal functions always holds locally and will hold globally under some causality condition \cite{bes1}. In the present situation we started from a set of observables that we did not know at first to be associated with some partial order. But we know now from (\ref{deforder}) that each observable of our set ${\cal C}$ is non-decreasing for $\preceq$. In the classical case the causal functions are exactly those which are non-decreasing for $\preceq$, but even there, we may very well start from a smaller set of causal functions which is still large enough to fulfill the hypothesis 1 and 2 and  end up with the same causal order. Hence we see that we need a ``saturation'' hypothesis to ensure that we started from the whole set of causal functions and not a smaller subset. The most natural is the following:

\begin{quote} {\bf Hypothesis 3  (saturation)} Every observable $o$ which have a non-decreasing Gelfand transform with respect to $\preceq$, i.e. $p_1\preceq p_2\Rightarrow p_1(o)\le p_2(o)$, is in ${\cal C}$.
\end{quote}

As above, the leap leading to this hypothesis cannot be fully justified. However we will give some arguments in the next section that we hope will give some credence to it.

\section{From the fundamental hypotheses to a set of working axioms}

In the previous section we put forward four hypotheses, which we qualify as fundamental, about the set ${\cal C}$ of causal observables. However, they are not very handy from a mathematical point of view. This is particularly true of the saturation hypothesis. In this section, we will see that under some very reasonable conjectures, they can be turned into axioms which are easier to work with.

First, we must specify the kind of Jordan algebra we use. In order to make sense of hypothesis 1 we must have a functional calculus. We choose JB-algebras, which possess continuous functional calculus. By doing so we subrepticely introduce a (noncommutative) topology, and the question of the compatibility with the causal structure arises. We will see in the next section that a very natural answer to this question can be given thanks to Gelfand theory for toposets.

Let us now use the saturation hypothesis. If $c_1,c_2\in {\cal C}$, then the Gelfand transforms $\hat c_1 : p\mapsto p(c_1)$ and $\hat c_2 : p\mapsto p(c_2)$ defined on pure states is non-decreasing, hence their sum is also. Therefore $c_1+c_2\in{\cal C}$ by the saturation hypothesis. Similarly, if $c_n$ is a sequence of elements of ${\cal C}$ converging in the norm to $c$, then each $\hat c_n$ is non-decreasing, hence the limit is, hence $c\in {\cal C}$.

Consider now the following set of axioms on ${\cal C}$, which can be derived from the fundamental hypotheses, as we have seen:

\begin{enumerate}
\item\label{ax1} ${\cal C}\subset {\cal O}$, a JB-algebra.
\item\label{ax2} ${\cal C}$ separates the pure states of ${\cal O}$.
\item\label{ax3} ${\cal C}$ is stable by non-decreasing continuous functional calculus.
\item\label{ax4} ${\cal C}$ is norm closed.
\item\label{ax5} ${\cal C}$ is stable by sum.
\end{enumerate}

Though it is a priori a weaker set of axioms than the fundamental hypotheses, we conjecture that it is in fact equivalent. We call this the \emph{saturation conjecture}. Whatever the status of this conjecture is, we will work in this paper with this maybe weaker set of axioms. In each case we will encounter, it is possible to show that the set ${\cal C}$ is indeed saturated, that is, contains all observables with non-decreasing Gelfand transform.

But first we want to strengthen axiom (\ref{ax2}) by asking instead that:

\begin{itemize}
\item [(ii')]${\cal C}$ separates the states of ${\cal O}$.
\end{itemize}

We will show below that axiom (ii) can be replaced with axiom (ii') without loss of generality provided the following property is true for ${\cal O}$:

\begin{quote}
{\bf Strong Stone-Weierstrass property} Let ${\cal O}$ be a JB-algebra and  $B$ be a subalgebra of ${\cal O}$. If $B$ separates $P({\cal O})\cup \{0\}$. Then $B={\cal O}$.
\end{quote}

The strong Stone-Weierstrass conjecture for JB-algebras states that all JB-algebras satisfy this property. This  is known to be true if $P({\cal O})$ is closed and also if ${\cal O}$ is post-liminal \cite{shep}. All algebras considered in this paper will be post-liminal.

\begin{propo} If ${\cal O}$ satisfies the Strong Stone-Weierstrass property, then the hypotheses (i) to (v) imply (ii'). 
\end{propo}
The proof is easy.  Let $B$ be the strong closure of the set of differences ${\cal C}-{\cal C}$. First notice that ${\cal C}-{\cal C}={\cal C}^+-{\cal C}^+$, where ${\cal C}^+$ denotes the set of positive elements in ${\cal C}$. Indeed, we can always rewrite a difference $a-b$ in the form $a+k.1-(b+k.1)$ where $k$ is a large enough real and $1$ is the unit of ${\cal O}$. (Remember that ${\cal C}$ contains the real multiples of $1$, as noticed in the previous section.) Next, consider $a,b\in{\cal C}^+$. Then $a+b\in{\cal C}^+$. Moreover, if $f$ is the continuous non-decreasing function such that $f(t)=0$ for $t<0$ and $f(t)=t^2$ for $t\ge 0$, then $f(a+b)=(a+b)^2$, $f(a)=a^2$ and $f(b)=b^2$ all belong to ${\cal C}^+$. Thus the Jordan product $a\circ b={1\over 2}((a+b)^2-a^2-b^2)\in B$. By linearity and continuity of the Jordan product we easily see that $B$ is indeed a subalgebra of ${\cal O}$. Now by hypothesis $B$ separates $P({\cal O})$, and since it contains a unit it separates $0$ from the other pure states. Thus $B={\cal O}$. Since $\overline{{\cal C}-{\cal C}}={\cal O}$ it is now obvious that ${\cal C}$ separates the states of ${\cal O}$.

\section{Isocones and $I^*$-algebras}

From now on we concentrate on JB-algebras which are the self-adjoint part of $C^*$-algebras. In this context the structure of the set of causal observables  becomes (under the saturation and Stone-Weierstrass conjectures) that of an isocone, the definition of which we recall below:

\begin{definition}  Let $A$ be a unital $C^*$-algebra, $\reel(A)$ be its self-adjoint part. An isocone of $A$ is a non-empty set $I$ such that:
\begin{enumerate}
\item $I\subset\reel(A)$,
\item\label{stabiso} $\forall \phi : \RR\rightarrow \RR$ continuous and non-decreasing, $ a\in I\Rightarrow \phi(a)\in I$,
\item\label{stabsom} $\forall a,b\in I$, $a+b\in I$,
\item\label{cl} $I$ is closed,
\item\label{sep} $I$ separates the states of $A$,
\end{enumerate}
If $I$ is an isocone of $A$, the couple $(I,A)$ is called an $I^*$-algebra.
\end{definition}

Isocones were first defined in \cite{bes1} but with slight differences, so let us make some comments on them to help make the connection.

First, instead of axiom (\ref{sep}) we required that the set of differences $I-I$ be dense in $\reel(A)$. This is seemingly a stronger condition, but in fact they are equivalent. Indeed, suppose $I$ separates ${\cal S}(A)$ and let $V=\overline{I-I}$. Note that $V$ is a real vector space by axioms (\ref{stabiso}) and (\ref{stabsom}). Suppose there exists $a\in\reel(A)$ such that $a\notin V$ and use Hahn-Banach theorem to find a continuous linear form $f$ such that $f(V)=0$ and $f(a)\not=0$. We can extend $f$ as an hermitian form on $A$ and by \cite{kr} there exist two states $s_1,s_2$ on $A$ and two nonnegative numbers $\alpha_1,\alpha_2$ such that $f=\alpha_1s_1-\alpha_2s_2$, and $\|f\|=\alpha_1+\alpha_2$. Since $1\in V$, we find that $\alpha_1=\alpha_2=\|f\|/2$. Hence $f=k(s_1-s_2)$ with $k=\|f\|/2$. By hypothesis $f(V)=0$, hence $s_1$ and $s_2$ are equal on $I$, and since $I$ separates the states, they are equal, which contradicts $f(a)\not=0$.

Moreover, in \cite{bes1} we did not put forward axiom (\ref{stabiso}), but instead asked that if $a,b\in I$ and $ab=ba$, then $\sup(a,b)$ and $\inf(a,b)$ must also be in $I$. This axiom was motivated by the hypotheses of the Kakutani-Stone theorem, thanks to which one can prove an analog of the Gelfand-Naimark theorem for commutative $I^*$-algebras (see below). However, this mathematically motivated axiom was later  shown to be equivalent to the more physically meaningful stability under non-decreasing continuous functional calculus in \cite{bes2}.

As alluded to above, the first step of the theory is to prove an analog of the Gelfand-Naimark theorem for commutative $I^*$-algebras. Let us introduce some terminology in order to state this theorem. Let $(M,\preceq)$ and $(N,\le)$ be two partially ordered topological spaces. A map $f : M\rightarrow N$ which preserves the order relations is said to be  \emph{isotone} or to be an \emph{isotony}. The set of continuous isotonies from $M$ to $\RR$, where $\RR$ is equipped with the natural order and topology, is denoted by $I(M)$. A \emph{toposet} is a topological ordered space which has ``enough'' isotonies, in the following sense:

\be
\forall x,y\in M, x\preceq y\Leftrightarrow (\forall f\in I(M), f(x)\le f(y))
\ee

This is the compatibility condition between  topology and order that we were wondering about in the previous section. A toposet is also called a completely separated topological ordered space. Note that if $M$ is compact, it is a toposet if and only if the relation $\preceq$ is closed as a subset of $M\times M$ \cite{nachbin}. A spacetime is a toposet if and only if it is causally simple \cite{bb}. It means that it does not contain closed causal curves and that the future and past of any event are closed. It is the second strongest causality condition after global hyperbolicity.

\begin{theorem}\label{GN} Let $M$ be a compact toposet. Then the couple $(I(M),{\cal C}(M))$ is a commutative $I^*$-algebra. Conversely, every commutative $I^*$-algebra $(I,A)$ is canonically isomorphic to one of this form, where $M$ is the set of pure states of $A$.
\end{theorem}

As we explained before, there cannot be any such duality theorem in the locally compact case where one uses the algebra ${\cal C}_0(M)$ of continuous functions vanishing at infinity since this does not contain any non-trivial isotone function in general. Hence one is forced to use a preferred compactification in this case.

The proof of theorem \ref{GN} can be found in \cite{bes1}. Let us just recall that the order on $M=P(A)$ is recovered from $I$ by the formula

\be
\forall \phi,\psi\in P(A), \phi\preceq_I\psi\Leftrightarrow (\forall a\in I, \phi(a)\le\psi(a))\label{deforder2}
\ee

which is just (\ref{deforder}) adapted to the present situation. As explained in the previous section, this formula  makes sense in any $I^*$-algebra, commutative or not, and we can see now that it endows $P(A)$ with the structure of a toposet. 
 
Before proceeding let us give a few examples of noncommutative $I^*$-algebras.

The first is trivial: if $I$ is the whole set of self-adjoint elements of $A$, all axioms are evidently satisfied. This is called the trivial isocone, and it induces on $P(A)$ the equality ordering $\phi\preceq \psi\Leftrightarrow \phi=\psi$.

For the second example we consider $A=M_2(\CC)$. In this case it is easily seen that any closed convex cone in $\reel A$ which contains the constants and has a non-empty interior is an isocone. We can give a geometric interpretation of the order induced by $I$ on $P(A)$ in the following way. First let us identify $P(A)$ with $2$-sphere of rank one projections which we call $S$. Then $K=I\cap S$ is compact, has non-empty interior and is either equal to $S$ or to a geodesically closed subset of a hemisphere. Then if we write $\omega_{[\xi]}$ for the pure vector state defined by the unit vector $\xi$ and $p_{[\xi]}$ for the projection on the line generated by $\xi$, we have

\be
\omega_{[\xi]}\preceq_I\omega_{[\eta]}\Leftrightarrow (\forall x\in K, d(x,p_{[\xi]})\ge d(x,p_{[\eta]}))
\ee

where $d$ is the geodesic distance on $S$. Clearly we can also use the identification of $P(A)$ with $\CC P^1$ and the Fubini-Study metric (see \cite{bes1} for details).

The next example is infinite-dimensional. If we consider any compact toposet $M$ with a regular Borel measure, we can faithfully represent the elements of ${\cal C}(M)$ as multiplicative operators on $L^2(M)$. Consider then the set $K$ of compact operators on $L^2(M)$ and $I=I(M)+\reel(K)$, that is the set of perturbations of isotonies of $M$ by self-adjoint compact operators. Then $I$ is an isocone of $A+K$ \cite{bes2}. In this case the pure state space of $A$ can be written $P(A)=X\coprod Y$ where $X\approx M$ and $Y$ is the set of vector states. Then the order induced by $I$ is the original order on $M$ extended trivially (no further relations on distinct elements) to $P(A)$. The intuitive idea behind this example is that compact operators correspond to infinitesimals, and that a non-decreasing function stays non-decreasing when one adds an infinitesimal to it. 

To build more examples, one can use the two following constructions.

\begin{theorem}\label{lexicosum} Let $(P,\preceq)$ be a \emph{finite} poset and for each $x\in P$ let $(I_x,A_x)$ be an $I^*$-algebra. We set $I=\bigoplus_{x\in P}I_x$, $A=\bigoplus_{x\in P}A_x$, and we write elements of $A$ in the form $(a_x)_{x\in P}$. We define 
$$\bigL_{x\in P}I_x=\{a\in I|\forall x,y\in P, x\prec y\Rightarrow \max\sigma(a_x)\le \min\sigma(a_y)\}$$
where $\sigma(a)$ denotes the spectrum of $a$. Then $\bigL_{x\in P}I_x$ is an isocone of $A$.
\end{theorem}

The isocone $L:=\bigL_{x\in P}I_x$ induces the lexicographic sum of the orderings $\preceq_x$ defined by $I_x$ on $P(A_x)$. More precisely put, we have $P(A)=\coprod P(A_x)$ and the pure states are ordered by $L$ in the following way: let $\phi\in P(A_x)$ and $\psi\in P(A_y)$ then
\be
\phi\preceq_L\psi\Leftrightarrow (x\not=y\mbox{ and }x\preceq y\mbox{ in }P)\mbox{ or }(x=y\mbox{ and }\phi\preceq_{I_x}\psi\mbox{ in }P(A_x))
\ee

This is why we call $L$ the lexicographic sum over $P$ of the isocones $(I_x)_{x\in P}$. If the order relation on $P$ is trivial, this is just a direct sum. There is a generalization of this construction to infinite posets, but one needs to add conditions on the order and the map $x\mapsto I_x$. To simplify matters we also suppose that the algebra $A_x=A$ is the same for all $x$. 

\begin{theorem}\label{lexicosuminfini} Let $(P,\preceq)$ be a metrizable compact poset and for each $x\in P$ let $I_x$ be an isocone of $A$. Let ${\cal A}$ be the $C^*$-algebra ${\cal C}(P,A)$ and 
$$\bigL_{x\in P}I_x:=\{f\in {\cal C}(A,P)|\forall x\in P, f(x)\in I_x\mbox{ and }\forall x,y\in P, x\prec y\Rightarrow \max\sigma(a_x)\le \min\sigma(a_y)\}$$
where $\sigma(a)$ denotes the spectrum of $a$. If $\prec$ is closed and $x\mapsto I_x$ is lower hemicontinuous, then $\bigL_{x\in P}I_x$ is an isocone of ${\cal A}$.
\end{theorem}

Remember that a multi-valued function $\varphi : X \longrightarrow 2^{Y}$ is said to be \emph{lower hemi-continuous}  if and only if the set $\{x \in X | \varphi(x) \cap V \ne \varnothing \}$ is open for all open $V \subset Y$. The condition that $\prec$ is closed as a subset of the product $P\times P$ is equivalent to the following one: for each $x\in P$ there exists a neighbourhood $U$ of $x$ such that no element of $U\setminus\{x\}$ is comparable to $x$. It is of course satisfied if $P$ is discrete (hence finite by compactness) but this is not necessary. For instance the ordering on $P=[0;1]$ defined by the relations $x\preceq 1$ for all $x\le 1/2$ and $x\preceq x$ for all $x\in[0;1]$ satisfies it. Note that if $\prec$ is closed and $P$ is compact then $\preceq$ is closed hence $P$ is a toposet.  

The other construction which will play an important role is the ``pushforward'' by surjective $*$-morphisms.

\begin{propo}\label{projinv} Let $A,B$ be $C^*$-algebras, $\pi : A\rightarrow B$ a surjective $*$-morphism, $I$ an isocone of $A$. Then $\overline{\pi(I)}$ is an isocone of $B$.
\end{propo}

Note that $\pi(I)$ is automatically closed if $B$ is finite-dimensional.

\section{Isocones in finite-dimensional and almost-commutative algebras}
The first task is to classify isocones in matrix algebras. We will say that an algebra is \emph{egalitarian} if it contains only the trivial isocone (the justification for the terminology is that the only isocone-induced order on its pure state space will be equality).

In \cite{bes2} we proved the following:

\begin{theorem} The matrix algebra $M_n(\CC)$ is egalitarian for all $n$ except $n=2$.
\end{theorem}

To put this result in perspective let us move away for a moment from the category of $C^*$-algebras.  First we can consider real $C^*$-algebras, which are important in the context of the spectral standard model.  Then there are two new cases of non-egalitarian algebras which are $M_2(\RR)$ and $M_2(\HH)$. However in the more natural arena of JB-algebras, the exceptional character of these cases disappear, as they all fit into the familly of spin factors. Remember that a spin factor, or Jordan algebra of Clifford type is the direct sum $N\oplus \RR$, where $N$ is a real Hilbert space, with Jordan product $(u,x)\circ (v,y)=(xu+yv,xy+\langle u,v\rangle)$. The Jordan algebras of self-adjoint real, complex and quaternionic $2\times 2$ matrices are isomorphic to spin factors with $N=\RR^2,\RR^3$ and $\RR^5$ respectively. It will be proven elsewhere that spin factors are the only simple finite-dimensional non-egalitarian JB-algebras  (with the possible exception of the Albert algebra, which is the only case still under investigation).

Returning to the case of $C^*$-algebras, we now give a classification result.

\begin{theorem}
Let $I$ be an isocone in a finite-dimensional $C^*$-algebra $A_F$. Then there exists a finite poset $P$ and positive integers $(n_x)_{x\in P}$ such that up to a $*$-isomorphism one has $A=\bigoplus_{x\in P}M_{n_x}(\CC)$ and $I=\bigL_{x\in P} I_x$. 
\end{theorem}

Since $M_{n_x}$ is egalitarian if $n_x\not=2$ and isocones in $M_2(\CC)$ are easily characterized, as we explained above, this theorem classifies all finite-dimensional $I^*$-algebras. It shows that the generalization of partial orders given by $I^*$-algebras is quite limited in the finite-dimensional case: partial order seems to be essentially a commutative phenomenon. However, the new freedom offered by the exceptional $n=2$ case could prove to be important in the context of the Standard Model, possibly in relation with T-symmetry breaking.

Let us make some comments on the physical interpretation of the above theorem. Recall that we see the elements of $I$ as observables corresponding to generalized clocks. If we admit the usual interpretation rule of quantum physics that eigenvalues of observables corresponds to values that can be possibly measured, then the order on $P$ can be characterized as follows: $x\preceq y$ if and only if no generalized clock ever give a result at $x$ which is strictly larger than its result at $y$. 

On a more mathematical stance, let us note that what the theorem tells us is that the partial order naturally defined on the pure state space by $I$ descends on the structure space, which is $P$. This was not guaranteed, and is proved by a combinatorial/analytical argument. If there is a conceptual reason for this, we are not aware of it.

We now turn to the important case of almost-commutative algebras, that is, $C^*$-algebras of the form ${\cal A}={\cal C}(M,A_F)$, where $A_F$ is a finite-dimensional algebra. We fix ${\cal A}$ an almost-commutative algebra with $M$ compact. For simplicity of exposition, we suppose that $A_F=M_n(\CC)$ with $n\ge 2$. The case of a general finite-dimensional algebra $A_F$ is not essentially different (we refer to \cite{bb} for details). Note that the pure state space of ${\cal A}$ is $M\times \PP(\CC^n)$.

\begin{theorem} Let $I$ be an isocone of $\cal A$. Then the order induced by $I$ on $M\times \PP(\CC^n)$ is lexicographic: there exists a partial ordering $\preceq_M$ on $M$ and for each $x\in M$ there exists an isocone $I_x$ of $M_n(\CC)$ such that : $(x,\xi)\preceq_I (y,\eta)\Leftrightarrow \left\{\matrix{x\prec_M y\cr \mbox{or}\cr  x=y \mbox{ and }\xi\preceq_{I_x}\eta}\right.$
\end{theorem}

The key to prove this theorem is the above classification of  finite-dimensional $I^*$-algebras and the pushforward construction of   section 4. Note that the ordering on $M$ is defined from $I$ through the following formula:

\be
x\preceq_I y\leftrightarrow \forall a\in I,\max(\sigma(a_x))\le \min(\sigma(a_y))\label{orderM}
\ee

Conversely, given a partial ordering $\preceq_M$ and a map $x\mapsto I_x$ into the isocones of $M_n(\CC)$, we can state necessary and sufficient conditions for the existence of an isocone of $\cal A$ inducing them.

\begin{theorem} If $M$ is metrizable and  $\preceq$ is  a lexicographic order on $M\times \PP(\CC^n)$ associated to a partial order $\preceq_M$ on $M$ and local isocones $I_x\subset M_n(\CC)$, then there exists an isocone $I$ in $\cal A$ inducing $\preceq$ iff
\begin{itemize}
\item $\prec_M$ is closed.
\item $L : x\mapsto I_x$ is lower hemi-continuous.
\end{itemize}
In that case, $I$ is the lexicographic sum of the isocones $I_x$ over $x\in M$.
\end{theorem}

This theorem would provide a classification of almost-commutative $I^*$-algebras (with $M$ metrizable) provided one could prove that the isocone inducing a given order is unique. This would of course follow from the saturation conjecture. Note that it is not difficult to prove that the lexicographic sum $I=\bigL_{x\in M}I_x$ is saturated, that is, it is exactly the set of self-adjoint elements of $\cal A$ whose Gelfand transform are isotone for $\preceq$. Hence, in the almost-commutative case the saturation conjecture amounts to say that   $I$ does not strictly contain any isocone inducing the same ordering on $M\times \PP(\CC^n)$.

Let us end this section by remarking that the closedness condition on $\prec_M$ is really imposed by the noncommutativity of $M_n(\CC)$ for $n\ge 2$. If $A_F$ is a more general finite-dimensional $C^*$-algebra then its structure space will have several sheets. If we consider $\CC\oplus M_2(\CC)$ for instance, there are two sheets  $M_1$ and $M_2$. Then, the order induced on   $M_1$ which is associated to the commutative algebra $\CC$ will not have the closedness restriction.

\section{Lorentzian spectral triple and causal cones}
Let us now summarize the approach initiated by Franco and Eckstein to the question of causality in noncommutative geometry. We refer to \cite{fe1} and \cite{fe2} for details, however we sometimes give indications of proofs which are  different than those given there, and might be in some cases easier to generalize.

One first defines Lorentzian spectral ``triples'', which in fact contains five pieces of information.

\begin{definition} A Lorentzian spectral triple is given by the following data: a pre-$C^*$-algebra without unit $A$, a unitization  $\tilde A$ of $A$, a Hilbert space $H$ on which  $A$ and $\tilde A$ are faithfully represented, an unbounded operator  $D$ with dense domain such that for all $a\in \tilde A$

\begin{itemize}
\item $[D,a]$ is bounded,
\item $a(1+\bra D\ket^2)^{-1/2}$ is compact, where $\bra D\ket^2={1\over 2}(DD^*+D^*D)$.
\end{itemize}
and $j$ is a bounded operator such that $j^2=1$, $j^*=j$, $[j,\tilde A]=0$, $D^*=-jDj$. Moreover one asks that  $j=-[D,T]$
 where $T$ is an unbounded operator such that $(1+T^2)^{-1/2}\in \tilde A$.
\end{definition}

The role of $j$ is to turn $H$ into a Krein space, and the requirement that $j$ can in fact be written $j=-[D,T]$ is the way Lorentzian signature is single out among arbitrary pseudo-riemannian ones.

The first example is the following: consider a complete and globally hyperbolic Lorentz manifold $M$, $S$ its spinor bundle, $H=L^2(S)$, $D=-i\gamma^\mu\nabla^S_\mu$ the usual Dirac operator. The algebras $A$ and $\tilde A$ are respectively ${\cal S}(M)$, the Schwartz functions on $M$, and the algebra generated by causal bounded functions with bounded derivatives. The role of the operator $T$ is played by a global time function on $M$, called $x^0$, so that $j=i\gamma(dx^0)$. On a coordinate chart where $x^0$ is the first coordinate one can write $j=i\gamma^0$. 

Given a lorentzian spectral triple, one can define a causal cone.

\begin{definition}
A causal cone ${\cal C}$ is a cone in the self-adjoint part of $\tilde A$ which contains $1$ and such that
\begin{itemize}
\item $\overline{span_\CC({\cal C})}=\overline{\tilde A}$
\item $\forall a\in {\cal C}$, $j[D,a]\le 0$ (is a negative semidefinite operator)
\end{itemize}
\end{definition}

We can then use a causal cone to define a partial order on the space of pure states $P(\tilde A)$, as one does with an isocone. We can see that we recover the usual partial order on $M$ (which is a subset of the compact space $P(\overline{\tilde A})$) through the following result (\cite{fe1}, theorem 11):


\begin{theorem} If $({\cal A},\tilde{\cal A},H,D,j)$ is the Lorentzian spectral triple built out of a globally hyperbolic Lorentzian manifold $M$, then $f\in\tilde{\cal A}$ is causal if and only if $j[D,f]\le 0$.
\end{theorem}

This can be most easily proved using a pseudo-orthonormal moving frame (vierbein). We can then define the Lorentz gradient of $f$ by $df(v)=g(v,\nabla f)$, and decompose $\nabla f$ in the form $-(\partial^0f)e_0+(\vec\nabla f)^ie_i$  in the moving frame.  We have  $i\gamma^0[D,f]=-\partial^0 f{\rm Id}+b$, with $b=\displaystyle{\sum_{1\le i\le 3}}\partial^i f\gamma^0\gamma^i$. Noting that $b^*=b$ and $b^2=\|\vec \nabla f\|^2{\rm Id}$, we see that the spectrum of $b$ is $\{\pm\|\vec\nabla f\|\}$. Hence $i\gamma^0[D,f]$ is negative semidefinite iff $\nabla f_x$ lies inside the past-cone at every $x\in M$, which is easily seen to be a necessary and sufficient condition for $f$ to be causal.

Now since globally hyperbolic manifolds are toposets, we see that the order defined by 
\be
x\le_M y\Leftrightarrow (\forall f\in\tilde{\cal A}, j[D,f]\le 0\Rightarrow f(x)\le f(y))
\ee
is exactly the original causal order of the manifold.

\section{Comparison of the two approaches}

We see that the two approaches agree as far as (globally hyperbolic) Lorentzian manifolds are concerned: a causal cone is just the smooth part of the isocone of causal functions. To see what happens in the noncommutative case, we will focus on a simple though enlightening example. We take $M=\RR^{1,1}$ to be the $2$-dimensional Minkowski spacetime. We consider the finite-dimensional spectral triple $(A_F,H_F,D_F)$ with $A_F=M_2(\CC)$, $H_F=\CC^2$, $D_F={\rm diag}(d_1,d_2)$. A product lorentzian spectral triple is defined by ${\cal A}={\cal S}(M)\otimes A_F$, ${\cal H}=L^2(S)\otimes H_F$, where $S$ is the spinor bundle which here is trivial with fiber $\CC^2$, $D=D_M\otimes 1+\gamma^0\gamma^1\otimes D_F$, with $D_M=-i\gamma^\mu\partial_\mu$, and $j=i\gamma^0\otimes 1$. The algebra $\tilde{\cal A}_M$ is defined to be the linear span of smooth causal bounded functions with bounded derivatives, and $\tilde{\cal A}=\tilde{\cal A}_M\otimes A_F$. This in fact corresponds to the Penrose compactification of the Minkowski plane which we briefly recall. Let $(x^0,x^1)$ be the canonical coordinates on $\RR^{1,1}$, and write $u=x^0+x^1$, $v=x^0-x^1$. The map $\phi : M\rightarrow P=]-\pi;\pi[^2$ defined by $\phi(u,v)=(\mu,\nu):=(2\tan^{-1}(u),2\tan^{-1}(v))$ is   a smooth diffeomorphism (the factor of $2$ is inserted to stick with usual conventions). Moreover $\phi$ sends the causal order $\le_M$ on $M$, which in the coordinates $(u,v)$ is just the product order on $\RR^2$, to the product order $\le_P$ on $P$. Since the derivatives of $\phi$ are all bounded, the push-forward  $\phi_* : f\mapsto f\circ\phi^{-1}$  sends $\tilde{\cal A}_M$ to the linear span of bounded causal (for $\le_P$) functions on $P$ with bounded derivatives. Such functions can be extended by continuity on $\bar P=[-\pi;\pi]^2$. To see it notice  that causal functions on $P$ can be expressed in the form $g(\mu)+h(\nu)$ with $g,h$ non-decreasing on $\RR$. Using the Stone-Weierstrass theorem we then immediately see that the linear span of bounded causal functions on $P$ with bounded derivatives are dense in ${\cal C}(\bar P)$. Hence the pure state space of $\tilde{\cal A}_M$ is homeomorphic to the Penrose compactification of $M$.

Let us look for the elements of the causal cone. A hermitian-matrix-valued function $\alpha$ is in ${\cal C}$ if $j[D,\alpha]=\gamma_0\gamma_\mu\otimes\partial^\mu \alpha-\gamma_1\otimes i[D_F,\alpha]\le 0$. With $\gamma_0=\pmatrix{0&i\cr i&0}$, $\gamma_1=\pmatrix{0&-i\cr i&0}$ this condition can be written in block-matrix form:

\be
\pmatrix{2\partial_u\alpha&[D_F,\alpha]\cr -[D_F,\alpha]&2\partial_v\alpha}\ge 0\label{freq}
\ee

In particular this condition implies that 

\be
\partial_u\alpha\ge 0, \partial_v\alpha\ge 0\label{derpos}
\ee

Integrating (\ref{derpos}) we see that if $x\le y$ (where $\le$ denotes the usual causal order on Minkowski space), then $\alpha(x)\le\alpha(y)$ as operators. Conversely, if $\alpha(x)\le\alpha(y)$ for all $\alpha\in{\cal C}$, then using $\alpha=$diag$(t_1,t_2)$, where $t_1,t_2$ are two time functions, one sees that $x\le y$. Hence we have proved that

\be
x\le_M y\Leftrightarrow \forall\alpha\in{\cal C}, \alpha(x)\le\alpha(y)\label{form11}
\ee

This formula can be put in contrast with  (\ref{orderM}) (where, we recall, $\preceq_M$ cannot be equal to $\le_M$). In particular we note that if $\lambda_1(x)\le\lambda_2(x)$ are the eigenvalues of $\alpha(x)$, then (\ref{derpos}) shows that $x\mapsto \lambda_i(x)$ are causal functions. It is however perfectly possible that $\lambda_2(x)>\lambda_1(y)$ with $x\le_M y$. If we keep the eigenvalue/observed value correspondance and interpret the elements of ${\cal C}$ as causal observables, we have here a generalized clock which can assign a later date to $x$ than to $y$, which is inconsistent with the observations of section 5. This shows that a causal cone is not an isocone in the noncommutative case. The mathematical origin of the difference between causal cones and isocones is that the former are not stable under non-decreasing functional calculus (this can be traced back to the fact that monotone functions are not necessarily operator-monotone). 

Finally we recall the main result of \cite{fe2}, which is the characterization of the order induced on the pure states of the form $(x,\omega_{[\xi]})$, with $x\in M$ and $\omega_{[\xi]}$ is the vector state $a\mapsto \bra \xi,a\xi\ket$ on $M_2(\CC)$. Let us fix some notations before stating the theorem: $d_{NC}(p,q)$ is Connes' noncommutative distance between two pure states on $M_2(\CC)$ induced by $D_F$, and $\ell(x,y)$ is the Lorentz distance between $x$ and $y$, that is to say the supremum of $\int_\gamma ds$ where $\gamma$ runs over all future-directed causal curves from $x$ to $y$ (it is set to $0$ if $y$ is not in the future of $x$). In the case at hand, $\ell(x,y)$ is the proper time elapsed from $x$ to $y$ measured by an inertial observer.

\begin{theorem}
$(x,\omega_{[\xi]})\preceq_{\cal C}(y,\omega_{[\eta]})$ iff $x\le_M y$ and  $\ell(x,y)\ge d_{NC}(\omega_{[\xi]},\omega_{[\eta]})$.
\end{theorem}

Note that the inequality implies that $d_{NC}(\omega_{[\xi]},\omega_{[\eta]})$ is finite. We recall that it is the case iff $\omega_{[\xi]}$ and $\omega_{[\eta]}$ are at the same latitude on the $2$-sphere, where the north and south poles are given by the eigenstates of $D_F$ (see \cite{PM}). Let us explain briefly why the condition is necessary. First we can suppose without loss of generality that $x=(0,0)$ and $y=(s,0)$. Then   $\ell(x,y)=s$. Using matrix-valued functions of the form $\alpha(t,r)=tI_2+A$, where $A$ is a constant hermitian matrix, we easily obtain from (\ref{freq}) that $\alpha\in{\cal C}$ iff $\|[D_F,A]\|\le 1$, which gives the result.

We now turn our attention to $I^*$-algebras. We work directly with the Penrose compactification, so the algebra is ${\cal A}={\cal C}(\bar P)\otimes M_2(\CC)$. The other datum we need is an isocone ${\cal I}$ of ${\cal A}$. Note the difference with the previous approach: the usual causal order was encapsulated in the choice of $D$ and $j$, wheras here it is induced by the isocone, and as such cannot be the usual causal order on $\bar P$, since the noncommutative nature of the algebra ${\cal A}$ imposes that the causal relations locally vanish. What we ask then, is that the isocone-induced order $\preceq_\Lambda$ on $\bar P$ depends on a constant $\Lambda$ and tends to the usual causal order when $\Lambda\rightarrow 0$. If we moreover ask  $\preceq_\Lambda$ to be Lorentz-covariant, it will be easily constructed.

We first define $\preceq_\Lambda$ on $\RR^{1,1}$ and then we will send it to $\bar P$. On $\RR^{1,1}$ we set:

\be
(x^0,x^1)\preceq_\Lambda ({x'}^0,{x'}^1)\Leftrightarrow x^0\le {x'}^0\mbox{ and }(x^0-{x'}^0)^2-(x^1-{x'}^1)^2\ge \Lambda^2\label{deforder3}
\ee

That is, we simply replace forward light-cones with forward hyperbolae of mass $\Lambda$. Changing to the coordinates $(u,v)$, and using the notations $\Delta u=u-u'$, $\Delta v=v-v'$, we easily see that this definition is equivalent to

\be
(u,v)\preceq_\Lambda (u',v')\Leftrightarrow \Delta u\ge 0, \Delta u\Delta v\ge \Lambda^2
\ee

Hence in $P$, using the same notation to denotes the partial order we obtain

\be
(\mu,\nu)\preceq_\Lambda (\mu',\nu')\Leftrightarrow \Delta \mu\ge 0, (\Delta\tan{\mu\over 2})(\Delta\tan{\nu\over 2})\ge \Lambda^2
\ee

with obvious notations, which is extended to $\bar P$ by $(\mu,\nu)\preceq_\Lambda (\mu',\nu')\Leftrightarrow \Delta\mu\ge 0,\Delta\nu\ge 0$ anytime $(\mu,\nu)$ or $(\mu',\nu')$ belongs to the boundary.

\begin{figure}[hbt]
\begin{center}
\includegraphics[scale=0.8]{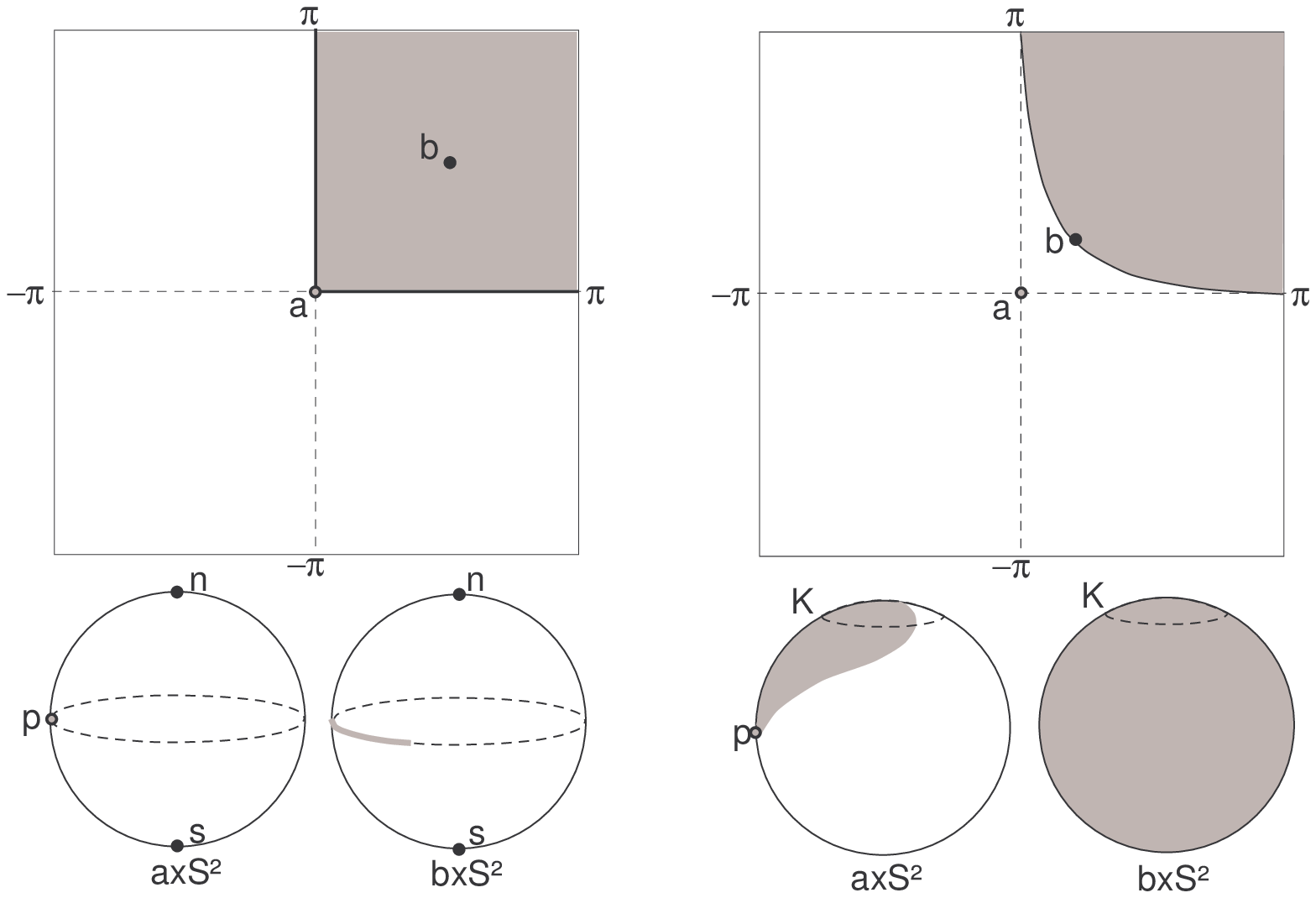}
\end{center}
\caption{This figure shows the set of pure states which are larger than $(a,p)$, where $a$ is an event in the Penrose compactification of the Minkowski plane, and $p$ is a pure state in $P(M_2(\CC))\simeq S^2$, for the two order relations described in the text: on the left side the order induced by a causal cone, and on the right from an isocone. Each point in the Penrose square represents a $2$-sphere. It is depicted in grey if the sphere contains at least a pure state larger than $(a,p)$. Under the square the two spheres at $a$ and $b$ are shown. On the left, only the states lying at the same latitude than $p$ can be in the future of $(a,p)$. The North and South poles are given by the eigenprojectors of the Dirac operator. In the sphere $b\times S^2$ one sees that the states in grey are those whose geodesic distance from $p$ is less that the Lorentz distance from $a$ to $b$. On the right we can see a non-trivial ``inner'' future light-cone in the sphere $a\times S^2$. It is easily seen to be the intersection of $S^2$ with $p+K^\circ$, where $K^\circ$ is the dual cone of $K$, which is here taken to be a spherical cap. At $b$ or any other point in the future of $a$ and not equal to $a$, the sphere is entirely in the future of $(a,p)$.}\label{fig1}
\end{figure}

To complete the data we have to choose a lower-hemicontinuous isocone map $(\mu,\nu)\mapsto I(\mu,\nu)\subset M_2(\CC)$. For instance we can take a constant map with value a non-trivial isocone in $M_2(\CC)$ defined by a closed  and geodesically convex subset $K$ of the sphere $S^2$. We finally define ${\cal I}=\bigL_{(\mu,\nu)\in \bar P}I(\mu,\nu)$, with $I_x=I=$constant, and $\bar P$ ordered by $\preceq_\Lambda$. We illustrate the two situations in figure \ref{fig1}.

\section{Conclusion}
It is evident from the above examples that there is much more freedom in the isocone approach, and it is not surprising since isocones generalize partial orders which can be arbitrary apart from possesing the toposet property. In particular we do not have any axiom which would ensure the compatibility with a given metric, incarnated in the noncommutative setting by the Dirac operator. This is clearly an important ingredient to add in the future development of the theory. On the other hand we believe that the approach is sufficiently general and physically well motivated to be trusted on its qualitative conclusions, in particular on the necessary vanishing of the causal relations at small scale (with no indication on the smallness scale, though). As a matter of fact, equivalent conclusions have been reached in \cite{klv} and \cite{asz} by studying the properties of the spectral action at high energy. The agreement between the isocone approach and the properties of the spectral action is all the more striking that the latter is formulated in the euclidean formalism. However, this action reproduces, under some cut-off energy scale, the one of the Standard Model which we know is correct on experimental ground. An optimistic explanation of the agreement between the disappearance of causality at small scale/high energy, as required by the isocone approach, and the non-propagation of high-energy bosons implied by the spectral action formalism could therefore be the following one: 1) the fundamental hypotheses laid up in the first section can be trusted (i.e. Nature really behaves like this), 2) the asymptotic developments of the spectral action can be trusted, 3) Nature is consistent with itself. To temperate this enthousiasm, let us observe that the local disappearance of causality implied by the isocone approach might be more soberly viewed as just another theorem of the noncommutativity-implies-discreteness kind  that we are familiar with since the advent of quantum theory.

However, as we have seen, causal cone theory does not predict such a discreteness. Quite the contrary, it is fully compatible with usual causality at all scale. Does it mean that it is wrong-headed ? Not necessarily. In fact it is interesting to note that the bottom-up isocone approach and the top-bottom causal cone theory each has what the other lacks. The first has, we believe, a clear physical interpretation and is able to derive sensible qualitative conclusions at large energies, the second has a direct contact with the usual apparatus of noncommutative geometry, providing in particular a very transparent connection between the Dirac operator and the causal order. We think that it is not only possible, but necessary, that the two approaches merge. It is all too obvious that the culprit of their present incompatibility is Minkowski spacetime, whose smoothness cannot be trusted at arbitrary large energy scale. In  replacing  smooth spacetime with a discrete structure we do not expect to encounter any difficulty on the isocone side, but we will have to discretize the Dirac operator along the lines of \cite{VSM}. Whether this permits to equate causal cones with isocones is under investigations.










\begin{thebibliography}{9}
\bibitem{cm} Chamseddine AH, Connes A, Marcolli M  2007, { Gravity and the Standard Model with neutrino mixing}, {\it  Adv. Theor. Math. Phys.} {\bf 11} 991 
\bibitem{vddg}   van den Dungen K,  Paschke M,  Rennie A 2013, { Pseudo-Riemannian spectral triples and the harmonic oscillator}, {\it J. Geom. Phys.} {\bf 73}  37 
\bibitem{besq} Besnard F 2007, {\it Canonical quantization and the spectral action, a nice example}, {\it J. Geom. Phys.} {\bf 57}  1757
\bibitem{dfr} Doplicher S, Fredenhagen K, Roberts  JE 1995, {The quantum structure of spacetime at the Planck scale and quantum fields}, {\it Commun. Math. Phys.}, {\bf 172} 187 
\bibitem{sorkin}  Bombelli L,  Lee J,  Meyer D,  Sorkin RD 1987, {  Space-time as a causal set}, {\it Phys. Rev. Lett.} {\bf 59}, 521 
\bibitem{bes1} Besnard F 2009, A noncommutative view on topology and order, {\it J. Geom. Phys.} {\bf 59} 861
\bibitem{fe1} Franco N, Eckstein M 2013, An algebraic formulation of causality for noncommutative geometry, {\it Class. Quant. Grav.} {\bf 30} 135007
\bibitem{tppm} Besnard F 2011, Time of philosophers, time of physicists, time of mathematicians, {\it Preprint} arXiv:1104.4551
\bibitem{tqm} Muga JG, Sala Mayato R, Egusquiza IL (Eds.) 2002, {\it Time in quantum mechanics},  (Berlin: Springer)
\bibitem{mingsan} Minguzzi E, S\`anchez M 2008, The causal hierarchy of spacetimes, {\it Recent developments in pseudo-Riemannian geometry}, ESI Lect. Math. Phys 299

\bibitem{emch} Emch GG 1984, {\it Mathematical  and Conceptual Foundations of 20-th Century Physics}, (Amsterdam: Elsevier)
\bibitem{ae} S. T. Ali, M. Engli\v{s}, {\it Quantization methods: a guide for physicists and analysts}, Rev. Math. Phys. {\bf 17} (2005) 391 math-ph/0405065
\bibitem{rov} Rovelli C 2002, GPS observables in General Relativity, {\it Phys. Rev. D} {\bf 65} 044017
\bibitem{shep} Sheppard B 2001, A Stone-Weierstrass theorem for postliminal JB-algebras.  {\it Q. J. Math.} {\bf 52} 507 

\bibitem{bes2} Besnard F 2013, Noncommutative ordered space: examples and counterexamples, {\it to appear in Class. Quantum Gravity}, arXiv:1312.2442
\bibitem{nachbin} Nachbin L 1965, {\it Topology and Order}, (Princeton: D. Van Nostrand Company)
\bibitem{bb} Bizi N, Besnard F 2014, The disappearance of causality at small scale in almost-commutative manifolds, {\it Preprint} arXiv:1411.0878

\bibitem{kr} Kadison R, Ringrose  J 1983, {\it Fundamentals of the Theory of Operator Algebras}, (New York: Academic Press)
\bibitem{fe2} Franco N, Eckstein M 2014, {Exploring the Causal Structures of Almost Commutative Geometries}, {\it SIGMA} {\bf 10}  010
\bibitem{PM} Iochum B, Krajewski T, Martinetti P 2001, Distances in Finite Spaces from Noncommutative Geometry, {\it J. Geom. Phys.} {\bf 37}, 100
\bibitem{klv} Kurkov MA, Lizzi F, Vassilevich D 2014,  High energy bosons do not propagate, {\it Phys. Lett. B} {\bf 731},   311	
\bibitem{asz} Alkofer N, Saueressig F, Zanusso O 2014, Spectral dimensions from the spectral action, {\it Preprint} arXiv: 1410.7999
\bibitem{VSM} Marcolli M, van Suijlekom WD 2014, Gauge networks in noncommutative geometry, {\it J. Geom. Phys.} {\bf 75}, 71

\end{thebibliography}
\end{document}